\documentclass{amsart}

\theoremstyle{definition}

\theoremstyle{remark}

\numberwithin{equation}{section}

\begin{document}


\title{Eye of the Beholder}


\author{*Nikolas Aksamit}
\address{Department of Mathematics, University of Utah, 
Salt Lake City, UT 84105}
\email{aksamit@math.utah.edu}
\address{Now at: Department of Geography and Planning, University of Saskatchewan, Saskatoon, SK S7N5C8, Canada}


\author{Don H. Tucker}



\begin{abstract}

We show that the real line $\mathbb{R}$ viewed as a vector space is of uncountable (algebraic) dimension over the scalar field $\mathbb{Q}$ of rational numbers. We then build an operator $J$ which maps $\{\mathbb{R}, \mathbb{Q}\}$ onto $\{\mathbb{R}, \mathbb{Q}\}$, is $\mathbb{Q}$-linear and whose graph is scattered all over the place, yet is still continuous in the inner product structures on the domain and range spaces. $J$ is not continuous if the usual norm is used on either the domain or range space. We lose continuity and linearity if the scalar field is completed.

\end{abstract}


 \maketitle





\section{Introduction}

Much of the mathematical structure humans have developed is predicated upon our intuitions which are very much planet earth based: counting our sheep, stepping off distances, using plumb bobs and crude levels to aid in construction of shelters, irrigation ditches, and so on. There is little doubt that such activities led to rational arithmetic and geometry. It should not be altogether surprising that such structures, even elaborated, enjoy a degree of effectiveness in addressing other issues viewed from similar vantage points by the eyes of the same beholders.

Every student of high school  algebra has heard of the Pythagoreans' discovery that the diagonal of the unit square does not have a rational number for its length. The student may not have been told that this fact shot a serious hole in the Pythagoreans' notions that the rational numbers were adequate for arithmetic and geometric needs. A major consequence of this discovery was that ruler and compass constructions do not work in such systems and our geometric, that is planet earth based, intuition is not as reliable as had been thought at that time. Natural philosophy suffered a major wound that was not satisfactorily healed until the reals were completed in the mid to late 1800s.

Ruler and compass construction is not the only problem involved with this issue. Much of what we teach our students about mathematical analysis (calculus sorts of things) is formulated in terms of coordinate systems, that is basis elements in a vector space over a scalar field. We usually tell the students that the space is of dimension the cardinality of the set of basis elements. Such statements are fraught with the risk of great misunderstanding. Most students who have had a class in linear algebra are aware that the complex numbers are of dimension two over the scalar field $\mathbb{R}$, but of dimension one over the scalar field $\mathbb{C}$. Hardly any would believe the geometric dimension of $\mathbb{C}$ is one, given the triangles and other sets of positive area in $\mathbb{C}$.

In as much as we frequently build metrics and topologies based upon the coordinates of vectors, our notions of closeness can vary greatly. It is not likely that one will find a text which makes an issue of these issues.

Our teaching methods tend to concentrate our attention on the vector objects, giving brief mention to the scalars which are a lens through which we view the matters of closeness.

In this note we illustrate that a given set of vectors may be of dimension one relative to one scalar field and of uncountable dimension relative to another, the set of vectors being the same in both cases. This leads to a wide disparity in the notions of closeness which arise naturally from the coordinates in these settings.

\section{Inner Product Structure}

Suppose each of $X$ and $Y$ is a vector space over the scalar field $\Phi$ and $L$ is a linear map from $\{X, \Phi\}$ into $\{Y, \Phi\}$. Denote by $L_0$, the kernel of $L$, that is all $x\in X$ such that $L(x)=\theta$ the identity for addition in $Y$; $L_0$ is a subspace of $X$.

By the span of a set $S$ in $X$ we mean all finite linear combinations of elements in $S$, that is all finite sums $\Sigma_{i=1}^n\phi_ix_i$ where $\phi_i\in\Phi$ and $x_i\in S$. The span of a set $S$ is a subspace of $X$ which we denote $span_\Phi S$

The statement that $\{x_\gamma\}_{\gamma\in\Gamma}$ is an algebraic (Hamel) basis for $X$ means that if $x\in X$, then 

$$
x\overset{!}{=}\sum_{i=1}^n\phi_i x_i
$$
where $\phi_i \in\Phi$, $x_i\in\{x_\gamma\}_{\gamma\in\Gamma}$, and the symbol $\overset{!}{=}$ means the representation is unique and moreover, the sum involves only finitely many $x_i\in\{x_\gamma\}_{\gamma\in\Gamma}$ even though $\Gamma$ may be uncountable.

We now give an explicit construction of an algebraic basis for $X$ and one for $Y$ which depend upon $L_0$ (though the same type of construction could be used on any vector space absent a map $L$).

Denote the origin (or additive identity) of $X$ by the symbol $\theta$. Use the axiom of choice to well-order $X$. This also well-orders $L_0$. Toss out $\theta$ and choose the first element remaining in the set $L_0$. Call it $x_{a_1}$. Call the first element in $L_0 \setminus span_{\Phi}\{x_{a_1}\}, x_{a_2}$ and give the name $x_{a_3}$ to the first remaining element in $L_0 \setminus span_{\Phi}\{x_{a_1}, x_{a_2}\}$ and continue. Suppose this process does not exhaust $L_0$. There is a first remaining element, add it to our list and do it again; so if any element in $L_0$ had been overlooked by our process, there would have been a first one and it was not overlooked. We have a linearly ordered algebraic basis for $L_0$. Consider the first element remaining in $X\setminus L_0$, denote it $\{x_{b_j}\}$. Continue the process to obtain a basis for $X\setminus L_0$, and then $\{x_{a_i}\} \cup \{x_{b_j}\}$ is an algebraic basis for all of $X$. If $x\in X$, then $x\overset{!}{=}\Sigma \alpha_i x_{a_i} + \Sigma \beta_j x_{b_j}$ where $x_{a_i}\in \{x_{a_i}\}$ and $x_{b_j}\in \{x_{b_j}\}$ and so $x$ has a unique decomposition $x\overset{!}{=} x_0 +x_0^c$ where $x_0\in L_0$ and $x_0^c\in span_\Phi\{x_{b_j}\}$. We denote $span_\Phi\{x_{b_j}\}$ by $L_0^c$ and refer to it as the complement of $L_0$, with apologies for the use of the word "the." Different well-orderings will give rise to different algebraic bases for $L_0$ and indeed different complements for $L_0$, but will not change $L_0$, which is determined by $L$.

Denote by $R_L$ the range of $L$ in $Y$. Since $L$ is linear, $R_L$ is a subspace of $Y$ and each non-$\theta$ element in $R_L$, say $y$, is $L(x)$ for some $x\in L_0^c$. Indeed, the set $\{L(x_{b_j})\}$ spans $R_L$ over $\Phi$ and is a basis if $L$ is injective. Otherwise, $\{L(x_{b_j})\}$ may not be a basis due to duplicity, in which case we can obtain a basis with the attributes we desire as follows.

Set $y_1=L(x_{b_1})$, toss out from $R_L$ all elements in $span_\Phi\{L(x_{b_1})\}$. If what remains is the empty set, then $\{L(x_{b_1})\}$ is a basis for $R_L$. If not, then one of the elements of $\{L(x_{b_j})_{j\neq1}\}$ remains. In fact, a first one remains. Choose it as our second basis element and continue the process. This exhausts $R_L$ because $R_L=L(L_0^c)$.

Note: for $y_R\in R_L, y_R\overset{!}{=}\sum_i\gamma_i y_i\overset{!}{=}\sum_i\gamma_iL(x_i)$. It's unique anyway because we did not use the duplications of $L(x_{b_j})$'s in selecting the $y_i$'s.

Now, extend this basis to all of $Y$ by well ordering $Y\setminus R_L$ and proceeding as before. It follows that if $y\in Y$, then $y\overset{!}{=}y_R+ y^c$ where $y_R\in R_L$ and $y^c$ is in the complement of $R_L$ as determined by the extended basis. This says that $y\overset{!}{=} \sum_i\gamma_iL(x_i) + \sum_j \varphi_jy_j$ where $y_j$ is in the extended basis and $L(x_i)$ belongs to our basis for $R_L$.

Denote by $X^*$ the set of all linear maps from $X$ into $\Phi$. It is called the algebraic dual of $X$. Also suppose $\{x_\gamma\}_{\gamma\in\Gamma}$, $\Gamma$ an indexing set, is an algebraic basis for $X$, any algebraic basis will do for now. Suppose $x^*\in X^*$ and $x\overset{!}{=}\sum^n\alpha_i x_i$. Then $x^*(x)=\sum^n\alpha_ix^*(x_i)$. Denote $x^*(x_i)$ by $\beta_i$ for each $x_i\in\{x_\gamma\}_{\gamma\in\Gamma}$. Thus $x^*$ corresponds to a unique $\Phi$-valued function on the index set $\Gamma$ and therefore the representation $\{\beta_i\}_{\i\in\Gamma}$ may well be uncountable. Even so, $x^*(x)=\sum^n\alpha_i\beta_i$ is a finite sum for each $x\in X$. Let's look further. For each $i\in\Gamma$, define $x_i^*(x_j)=\delta_{ij}$ where $\delta_{ij}=0$ if $i\neq j$ and $\delta_{ii}=1$ for $i, j\in\Gamma$. Now suppose $x\overset{!}{=}\sum\alpha_i x_i$, then $x_j^*(x)=\sum\alpha_i x_j^*(x_i)=\alpha_j$ and therefore $\alpha_i=x_i^*(x)$. The coefficients of $\alpha_i$ are usually called the Fourier coefficients of the vector $x$ relative to the set $\{x_i^*\}_{i\in\Gamma}$ in $X^*$. We shall refer to the set $\{x_i^*\}_{i\in\Gamma}$ as a functional basis for $X^*$. In general, this set will not be an algebraic basis for $X^*$.

The two structures we have constructed on $X$ and $X^*$, namely the algebraic basis on $X$ and the functional basis on $X^*$ serves as mechanisms to define an inner product on $X$. Assume the elements of $\Phi$ are for example only either in $\mathbb{Q}, \mathbb{R}$ or $\mathbb{C}$, $x\overset{!}{=}\sum^n\alpha_ix_i$ and $y\overset{!}{=}\sum^m\beta_jx_j$. Define $\langle x, y \rangle = \langle \sum^n\alpha_ix_i, \sum^m\beta_j x_j\rangle=\sum^m\overline{\beta_j}x_j^*\sum^n\alpha_ix_i$ which after expansion becomes 

$$
\langle x, y \rangle = \sum_{i, j}\overline{\beta_j}\alpha_jx_k^*(x_i)=\sum\overline{\beta_j}\alpha_j
$$
and may well be zero, especially if no $x_i=x_j$ in the expansions of $x$ and $y$.

The symbol $\langle \quad ,\quad\rangle$ is linear in the first position, conjugate linear in the second if $\Phi=\mathbb{C}$, and bilinear if $\Phi=\mathbb{R}$. It easily follows that, $\langle \quad ,\quad\rangle$ is an inner product on $X$.




We define a norm on $X$ in the usual way: $\|x\|_X=\langle x,x\rangle^{\frac{1}{2}}=[\sum^n|\alpha_i|^2]^{\frac{1}{2}}$ where the $\alpha_i$ are the Fourier coefficients of $x$, i.e., $x\overset{!}{=}\sum\alpha_ix_i$. It follows that $\{x_\gamma\}_{\gamma\in\Gamma}$ is an orthonormal set.

This same type of construction on $Y$ with the algebraic basis we have developed (using $L$) nets us $Y^*$ plus an inner product and a norm on $Y$. That is if $y\overset{!}{=}\sum_i\gamma_iL(x_i) + \sum_j \varphi_j y_j$, then $\|y\|_Y=\langle y, y \rangle^{\frac{1}{2}} = [\sum_i|\gamma_i|^2 + \sum_j|\varphi_j|^2]^{\frac{1}{2}}$

It follows that the specific $\{L(x_i)\}$ we selected and $\{y_j\}$ form an orthonormal basis for $Y$. With these inner product induced norms on $X$ and $Y$, if $L$ is injective, it is a bounded linear operator of norm 1. (However, injectivity is not a necessary condition.)

To see why, consider for $x\in X, x= x_0+x^c$ where $x_0\in L_0, x^c\in L_0^c, x^c \overset{!}{=} \sum_{j=1}^m\beta_jx_j$ and  $x_0\overset{!}{=}\sum_{k=1}^n\alpha_kx_k$. We have $L(x)=L(x_0+x^c)=0+L(x^c)=\sum_{j=1}^m\beta_jL(x_j)$. $\|x\|_X=(\sum|\alpha_k|^2 + \sum|\beta_j|^2)^{\frac{1}{2}}$ because $L$ is injective. When building our basis for $Y$, each $L(x_i)$ is a basis element for $R_L$ and none are thrown out, therefore, $\|L(x)\|_Y=\|\sum_{j=1}^m\beta_jL(x_j)\|_Y=(\sum|\beta_j|^2)^{\frac{1}{2}}\le(\sum|\alpha_k|^2+\sum|\beta_j|^2)^{\frac{1}{2}}$. So, $\|L(x)\|_Y\le\|x\|_X$ for every $x\in X$ with equality if $x\in L_0$.

The observation that injectivity is not necessary will come in  the form of an unsettling example.


\section{The Example}

We will construct a linear map $J$ from a vector space $X$ to a vector space $Y$,  $X$ and $Y$ each being $\{\mathbb{R}, \mathbb{Q}\}$. The construction of $J$ involves building an algebraic basis for $X$. We will then build bases and inner product structures on $X$ and $Y$ induced by $J$ as outlined before. (The basis for inner product on $X$ is not necessarily the same basis we defined for $J$.) Lastly, we explore how the topological perspective of $J$ changes using our inner product based norm topology, the absolute value topology, or both on $X$ and $Y$.


\subsection{Example}

Suppose $X=\{\mathbb{R}, \mathbb{Q}\}$ and $Y=\{\mathbb{R}, \mathbb{Q}\}$. Well order $X=\mathbb{R}$ and construct an algebraic basis as before. Since there are uncountably many irrationals in $\mathbb{R}$ and our scalar field $\mathbb{Q}$ is countable it will require uncountably many elements to construct our basis. Call it $\{h_\gamma\}_{\gamma\in\Gamma}$. Choose any countable subset $H$ of this basis. Our subset $H$ is linearly ordered. Therefore we can label them $\{h_i'\}_{i\in\mathbb{N}}$.

Rescale $h_2'$ by a rational number $q$ such that $|qh_2'|<|\frac{h_1'}{2}|$. Rename $h_1'=j_1$ and $qh_2=j_2$. Replace $h_2'$ with $j_2$ in our basis. Note that $qh_2'\in span_{\mathbb{Q}}\{j_2\}$ implies $span_{\mathbb{Q}}\{h_2'\}=span_{\mathbb{Q}}\{j_2\}$ and we still have an algebraic basis. Continue this process replacing basis elements in such a way that $|j_{n+1}|<|\frac{j_n}{2}|$ for every $n\in\mathbb{N}$. Then $\{j_i\}_{i\in\mathbb{N}}\cup\{h_\gamma\}_{\Gamma\setminus\mathbb{N}}$ is an algebraic basis for $X$.

Define $J: X\to Y$ to be the sum of the (rational) Fourier coefficients in this unique representation. If, $x\overset{!}{=}\sum_{j=1}^mq_jx_j$ then $J(x)=\sum_{j=1}^mq_j$. It follows that $J$ is $\mathbb{Q}$-linear. 

To see why, consider: If $x=\sum_iq_ix_i$ and $y=\sum_jq_jx_j$ with $\alpha\in\mathbb{Q}$. 

$$J(x+\alpha y)=J(\sum_iq_ix_i+\alpha\sum_jq_jx_j)$$
$$=J(\sum_{i, j}q_ix_i +\alpha q_jx_j)=\sum_iq_i +\alpha\sum_jq_j=J(x)+\alpha J(y)$$

We denote by $J_0$ the real numbers whose coefficients (in the basis that we constructed for $J$) sum to zero. Construct an algebraic basis for $J_0$ as initially outlined by well ordering. Call this basis $\{k_\alpha\}$. Consider the remaining elements in $X\setminus J_0$. The element $j_1$ is surely in $X\setminus J_0$ because $J(j_1)=1$, and indeed $span_\mathbb{Q}\{j_1, J_0\}=X$, that is, $J_0$ is of codimension 1.

Let's check. Suppose there exists an element in $X\setminus span_\mathbb{Q}\{j_1, J_0\}$. Call it $a$, then $J(a)=p$ for some $p\in\mathbb{Q}-\{0\}$. Set $b=-pj_1 + a$; $b\in J_0$ because $J(b)=J(-pj_1+a)=J(-pj_1)+J(a)=-p+p=0$. Therefore, $a=b+pj_1$ and $a\in span_\mathbb{Q}\{j_1, J_0\}$, a contradiction.

It follows that if $x\in X$, then $x=qj_1+x_0$ where $qj_1\in J_0^c$, $x_0\in J_0$. Moreover if $x\in X$ and $J(x)=p$, then $x\overset{!}{=}x_0+p j_1$ because $$p=J(x)=J(qj_1+x_0)=J(qj_1)+J(x_0)=qJ(j_1)=q.$$

Notice that $R_J=J(X)=\mathbb{Q}\subset Y$. We do not yet have a basis for $Y$. Let's build one. Well-order $Y$ and build a basis for $Y$ by choosing $J(j_1)$ to be our first basis element. $span_\mathbb{Q}\{J(j_1)\}=\mathbb{Q}=R_J=span_\mathbb{Q}\{1\}$. Define the rest of our basis for $Y$ as we have done before, choosing first elements in $\mathbb{Q}$-span complements. We see that our basis for $Y$ is $\{J(j_1), y_\gamma\}_{\gamma\in\Gamma}$. We now have that $J:X\to Y$ is continuous with the outlined inner product structures on $X$ and $Y$.

To see why, recall that for linear operations between normed linear spaces boundedness is equivalent to continuity. Suppose

$$x\overset{!}{=}qj_1+\sum_{i=1}^nq_ik_{\alpha_i}\quad \|x\|_X=(|q|^2+\sum|q_i|^2)^{\frac{1}{2}}$$ then

$$J(x)=J(qj_1)+J(\sum q_ik_{\alpha_i})=qJ(j_1)+0=q$$ 

Thus, $\|J(x)\|_Y=|q|\le(|q|^2+\sum|q_i|^2)^{\frac{1}{2}}$ and $J(x)$ is bounded (by 1) and therefore continuous. However, $J: X\to Y$ is not a continuous operator with the absolute value topologies on $X$ and $Y$.

To see this, let $\{z_n\}$ be a sequence in $X$ such that $z_i=j_i$ for every $i\in\mathbb{N}$ ($\{j_i\}$ being the countable subbasis we constructed for $J$) As a sequence, $z_n=o(n).$ In fact, $z_n=O(\frac{1}{2^n})$. However, $J(z_n)=J(1\cdot j_n)=1\cdot J(j_n)=1$ for every $n\in\mathbb{N}$. Thus, $\{z_n\}$ is a sequence in $X$ that converges to $0$, but $J(z_n)$ is constantly one, hence converges to $1$ in $Y$. Therefore, $J(z_n) \not\rightarrow J(0)$.

If $X$ is given the inner product topology we constructed and $Y$ the absolute value topology, then $J:X\to Y$ is not bounded, hence not continuous. To see why, suppose there exists a bound $a\in\mathbb{R}$ for $J$. That is $|J(x)|\le a\|x\|_X$ for every $x\in X$. If $x=j_1 + \cdots + j_n$ then $|J(x)|=n\le a\sqrt{\sum^n|1|^2}=a\sqrt{n}$. This implies $n^2<a^2n$ and $n<a^2$ for every $n\in\mathbb{N}$. We arrive at our contradiction and $J$ is not bounded, hence not continuous.

If $X$ is given the absolute value topology and $Y$ the inner product topology we constructed, $J:X\to Y$ is not a continuous operator. This follows from the following: let $\{z_n\}$ be the sequence in $X$ as in the fourth paragraph back, then $\{z_n\}\rightarrow0$ in $X$. $J(z_n)=1$ in $Y$ for every $n\in\mathbb{N}$. $\|J(z_n)\|_Y=\|1\|_Y$ for every $n\in\mathbb{N}$. $\|1\|_Y\neq0$ because $\|\quad\|_Y$ is a norm. The result follows as before: $J(z_n)\not\rightarrow J(0)$.

Upshot: $J$ is a $\mathbb{Q}$-linear operator whose range is "scattered all over the place." J is  continuous in the inner product structure on $X$ and $Y$ and yet nowhere continuous if the absolute value topology is present on either $X$ or $Y$ (or both of course), and the other has the inner product topology.

It is worth noting if we restrict the domain of $J$ to be $\mathbb{Q}$, then $J$ is bounded and continuous with any of the four combinations of topologies by the following:

For any $q\in\mathbb{Q}$, $q=aj_1$ and $J(q)=J(a j_1)=aJ(j_1)=a$. Then

$$\|q\|_X=|a|=|J(q)|=\|J(q)\|_Y$$
$$|q|= b|J(q)|=b\|J(q)\|_Y$$

where $b=|j_1|$

It is not unusual that restricting the domain of a function so as to exclude points at which it fails to have a particular property can sometimes have happy circumstances.  What may have an unintended consequence is if we complete the scalar field $\mathbb{Q}$, but do not change the function, in our case $J$, then desirable properties may be lost, $J$ loses its linearity. While $J(\pi)$ still makes sense, it is not $\pi J(1)$.

\section{Conclusion}

When we construct a mathematical model of events that we observe it is our wont to identify certain parameters which we assume influence, perhaps effect, the events observed. An application of the Buckingham Pi Theorem is sometimes employed to reduce the set of parameters to a linearly independent subset which we assume will span the entire space of parameter values. We then assume a coordinate system in which this set of linearly independent parameters is identified as the basis elements. One challenge to the modeler and the model tester is to choose parameters which we can physically measure. We then assume the scalar field is $\mathbb{R}$ or $\mathbb{C}$ in order to achieve a topologically complete finite (algebraic) dimensional space in which we can (and do) attack the problem using the analysis we have developed for such situations.

The heart of the issue is the following. Our measurements are rational numbers and our computations with them yield only rational numbers. If our scalars are rational, then the space spanned by our finite set of parameters with those scalars is wanting in attributes needed to justify the applications of analysis. For example, if we have just one parameter, and it is discrete, and the scalars are $\mathbb{Q}$, then our space is of dimension 1, but we lack completeness. If the set of parameter-values is assumed to be a continuum and the scalars are $\mathbb{Q}$, we have a space of uncountable dimension. If the parameter values set and the the scalar field are completed to $\mathbb{R}$, we again have dimension one and analysis is available. However, as seen above, this may result in the loss of other desirable attributes such as linearity.


\begin{thebibliography}{20}

\bibitem{Naylor}  Arch W. Naylor and George R. Sell, \textit{Linear Operator Theory in Engineering and Mathematical Sciences},
Springer, 2000.

\bibitem{Sierpinski}  Waclaw Sierpinski, \textit{Cardinal and Ordinal Numbers},
Polish Scientific Publishers, 1965.

\bibitem{von Neuman}  John von Neumann, The Mathematician \textit{The World of Mathematics} {\bf 4},
Simon and Schuster, 1956.

\bibitem{Wigner} Eugene Wigner, The Unreasonable Effectiveness of Mathematics in the Natural Sciences,
\textit{Communications in Pure and Applied Mathematics} {\bf 13}, No. I, John Wiley and Sons,  1960.

\end{thebibliography}
\end{document}